\documentclass[12pt]{article}
\usepackage{latexsym}
\usepackage{amsfonts}
\usepackage{amsmath}
\usepackage{amssymb}
\usepackage{amscd}
\usepackage{bbm}

\newcommand{\diam}{\mathop{\rm diam}}

\newcommand{\R}{\mathbb{R}}

\newcommand{\inr}[1]{\left<#1\right>}

\newcommand{\E}{\mathbb{E}}

\newcommand{\eps}{\varepsilon}
\newcommand{\ep}{\varepsilon}

\newcommand{\conv}{\mathop{\rm conv}}
\newcommand{\supp}{\mathop{\rm supp}}

\newcommand{\Pp}{\mathbb{P}}
\newcommand{\uup}{uup}

\newcommand{\lstar}{\ell_*}

\makeatletter
\@addtoreset{equation}{section}

\newtheorem{Theorem}{Theorem}[section]
\newtheorem{Lemma}[Theorem]{Lemma}

\newtheorem{Definition}[Theorem]{Definition}

\newtheorem{Corollary}[Theorem]{Corollary}
\newtheorem{Remark}[Theorem]{Remark}

\newtheorem{Question}[Theorem]{Question}
\numberwithin{equation}{section} 

\def \proof {\noindent {\bf Proof.}\ \ }

\def \endproof
{{\mbox{}\nolinebreak\hfill\rule{2mm}{2mm}\par\medbreak}}

\date{August 21, 2006}
\title{Uniform uncertainty principle for Bernoulli and subgaussian
  ensembles}
 \author{
Shahar MENDELSON${}^{1}$ \and Alain PAJOR$^{2}$
       \and
Nicole TOMCZAK-JAEGERMANN${}^{3}$ }

\begin{document}
\maketitle \footnotetext{{\it 2000 MSC-classification:}
    46B07,  41A45, 94B75, 52B05, 62G99}
\footnotetext[1]{Partially supported by an Australian Research
Council Discovery grant.}
 \footnotetext[2]{Partially supported by
an Australian Research Council grant.}
\footnotetext[3]{This author holds the Canada Research Chair in
Geometric Analysis.}


\section{Introduction}
\label{sec:intro}


In \cite{CT1} Candes and Tao studied problems of approximate and exact
reconstruction of sparse signals from incomplete random measurements
and related them to the eigenvalue behavior of submatrices of matrices
of random measurements. In particular they introduced the notion they
called the {\em uniform uncertainty principle} (UUP, defined below)
and studied it for Gaussian, Bernoulli and Fourier ensembles.
This notion was further refined in \cite{CT2, CRT}.
In this context they asked
(\cite{T})
whether rectangular $k \times n$ Bernoulli
matrices (with $k < n$) have the property that by arbitrarily
extracting $m$ (with $m < k$) columns one can make
so obtained submatrices  arbitrarily close to (multiples of)
isometries  of a Euclidean space (of course $m$ would then depend on
the required degree of ``closeness'' and dimensions $k$ and $ n$).

A different--geometric--approach to approximate and exact
reconstruction problems was proposed in \cite{MPT1, MPT2}. Although
in these articles the notion of UUP was not considered, an
application of one of the main general results there in a simple
particular case implied an immediate affirmative answer to the
Candes-Romberg-Tao's
question (see Corollary 3.5 in \cite{MPT2} and
the comments afterwards).

\medskip

The common roots of the geometric approach of \cite{MPT1, MPT2}, as
well as
the UUP or other related properties, revolve around the fact that
various ``random projection'' operators may act as ``almost norm
preserving'' on various subsets of the sphere; with the UUP
associated to the  subset  of ``sparse'' vectors on the sphere
(denoted later by $U_m$).

In this note we observe that the results of \cite{MPT1, MPT2} can be
applied to a number of other sets with a very simple geometry to get
interesting conclusions for the Gaussian, Bernoulli, and more
generally, any subgaussian ensemble.
Since the proofs of the general results of \cite{MPT1, MPT2} are not
easily accessible to non-specialists, we also provide an alternative
elementary argument, which works for the specific sets we are
interested in.

\bigskip

Let us recall the following notation.  By $|\cdot|$ we denote the
Euclidean norm on $\R^n$, by $\inr{\cdot,
  \cdot}$ the corresponding inner product, by $B_2^n$ and $S^{n-1}$
the unit Euclidean ball and the unit sphere, respectively.  For $x=
(x_i)_{i=1}^n \in \R^n$ we let $\supp x = \{i: x_i \not = 0 \}$. For
a finite set $A$, the cardinality of $A$ is denoted by $|A|$, and
for a set $A \subset \R^n$, $\conv A$ denotes the convex hull of
$A$. Throughout, all absolute constants are fixed, positive numbers,
which are denoted by $c$, $C$, $c'$, etc. Their value may change
from line to line.

We will work with the following (slightly refined) definition of the
uniform uncertainty principle (\cite{CT1}).

\begin{Definition}
  \label{uup}
  A  $k \times n$ (random) measurement matrix $\Gamma$ obeys the uniform
  uncertainty principle with accuracy $0<\theta <1$ and oversampling factor
  $\lambda >1$, if the following statement is true with probability
  close to 1: for all subsets $A \subset \{1, \ldots, n\}$ with $|A|
  \le k/\lambda$, the matrix $\Gamma_A$, obtained
by extracting from $\Gamma$ the columns corresponding to $A$,  satisfies
\begin{equation}
  \label{eigenvalues}
  1 - \theta   \le {\lambda_{\min}\left(\frac{\Gamma_A^*\Gamma_A}{k}\right)}
  \le {\lambda_{\max}\left(\frac{\Gamma_A^*\Gamma_A}{k}\right)}\le 1+\theta,
\end{equation}
where $\lambda_{\min}$ and $\lambda_{\max}$ denote the minimal and maximal
eigenvalues, respectively.
Equivalently,
\begin{equation}
  \label{vectors}
  (1 - \theta)  |x|^2 \le \frac{|\Gamma x|^2}{k}
\le (1+\theta) |x|^2,
\end{equation}
for all vectors $x \in \R^n$ with $|\supp x| \le k/\lambda$.
\end{Definition}
We shall use a shorthand notation of $\uup(\theta, \lambda)$ for
the above property.

In this language and for the Bernoulli ensemble, Candes and Tao showed
(\cite{CT1, CT2}) that there exist two absolute constants $0 <
\theta_0 <1$ and $c >0$ such that for all $k < n$, $k \times n$
Bernoulli random matrices satisfy $\uup(\theta_0, \lambda)$ for
$\lambda=c\log(cn/k)$, and they asked (\cite{T}) whether an analogous
result is true for {\em every} $0 < \theta < 1$.  We formally state
their question as follows:

\begin{Question}
  \label{qu:CRT}
  Let
  $1 \le k <n$ and set $\Gamma$ to be a $k \times n$ random Bernoulli
  matrix. Let $0<\theta<1$ be arbitrary. Can one find $\lambda $
  depending on $\theta, k$ and $n$ only and satisfying $1 \leq \lambda
  \leq c_\theta\log(c_\theta n/k)$, where $c_\theta >0$ depends only on
  $\theta$, such that for $n$ ``large enough" and any $k < n$,
  $\Gamma$ satisfies $uup(\theta,\lambda)$ with probability close to
  1?
\end{Question}

\medskip

As already mentioned earlier, a positive answer to this question for
(more general) subgaussian measurements follows immediately from the
results of \cite{MPT1,MPT2}, and we explain this connection in the
next section. We then show how one can obtain a similar estimate using
elementary methods, which can also be used to solve the approximate
reconstruction problem in certain simple (but central for the
applications) cases (see Section \ref{oper_concentr} for more
details).

\bigskip

\noindent {\em Acknowledgement:\ } Part of the work on this article was conducted
during the Trimester ``Phenomena in High Dimensions'' held at the
Centre Emile Borel, Institute Henri Poincar\'{e}, Paris in the Spring
of 2006. We are grateful to the Institute for its hospitality
and the excellent working atmosphere it provided.

\section{Subgaussian matrices and geometry of the set of sparse vectors}
\label{subgaussian}

We first recall a few definitions.
Let $X$ be a random vector in $\R^n$; $X$ is called isotropic if
for every $y \in \R^n$, $\E |\inr{X,y}|^2 = |y|^2$, and is
$\psi_2$ with a constant $\alpha$ if for every $y \in \R^n$,
$$
\|\inr{X,y}\|_{\psi_2} := \inf\left\{s : \E \exp\left(\inr{X,y}^2
/s^2\right)  \le 2 \right\} \leq \alpha |y|.
$$
The most important examples for us are the Gaussian vector $(g_1,
\ldots, g_n)$ where the $g_i$'s are independent $N(0,1)$ Gaussian
variables and the random sign vector $(\ep_1, \ldots, \ep_n)$ where
the $\ep_i$'s are independent, symmetric $\pm 1$ (Bernoulli) random
variables; in both these cases the random vectors are isotropic with a
$\psi_2$ constant $\alpha =c_0'$, for a suitable absolute constant
$c_0' \geq 1$.

A subgaussian or $\psi_2$ operator is a random operator $\Gamma: \R^n
\to \R^k$ of the form
\begin{equation}
   \label{Gamma}
   \Gamma=\sum_{i=1}^k \inr{X_i,\cdot }e_i,
\end{equation}
where $X_1, \ldots, X_k$ are independent copies of an isotropic
$\psi_2$ vector $X$ on $\R^n$.

Note that if $X_i=(x_{i,j})_{j=1}^n$ then $\Gamma$ is represented by
a matrix whose rows are $(X_i)_{i=1}^k$. However, although the rows
of the matrix are independent random vectors, the entries within
each row may be dependent.

Finally, for a subset $T \subset \R^n$ we set
\begin{equation}
  \label{lstar_def}
\lstar(T)=\E \sup_{t
  \in T} \left| \sum_{i=1}^n g_i t_i\right|,
\end{equation}
where $ t = (t_i)_{i=1}^n \in \R^n$ and $g_1,...,g_n$ are independent
$N(0, 1)$ Gaussian random variables.

\medskip

The following fact was proved in \cite{MPT2} (Corollary 2.7) as a
consequence of one of the main results of \cite{MPT1, MPT2}.

\begin{Theorem}
  \label{expect}
Let $1 \le k \le n$ and $0<\theta <1$.
  Let $X$ be an isotropic $\psi_2$ random vector on $\R^n$ with
  constant $\alpha$, set $X_1, \ldots, X_k$ to be independent copies of
  $X$, put $\Gamma$ as defined by (\ref{Gamma}) and let
  $T \subset S^{n-1}$. If $k$
  satisfies
  $$
  k \ge (c' \, \alpha^4/\theta^2) \lstar (T )^2,
  $$
  then with probability at least $1- \exp(- \bar{c}\, \theta^2 k/\alpha ^4)$,
  for all $x \in T$,
\begin{equation}
  \label{two_sided-2}
   1 - \theta \le \frac{|\Gamma x|^2}{  k} \le 1+\theta,
\end{equation}
where $c', \bar{c} >0$ are absolute constants.
\end{Theorem}

Let us explain the meaning of Theorem \ref{expect}, and for the sake
of simplicity, assume that $\alpha$ is an absolute constant (in
particular independent on the dimension $n$), as this is the situation
for Gaussian or Bernoulli random vectors. The parameter $\lstar(T)$ is
a complexity measure of the set $T$; in this context, it measures the
extent in which probabilistic bounds on the concentration of
individual random variables of the form $|\Gamma x|^2$ around their
mean can be combined to form a bound that holds uniformly for every $x
\in T$. The assertion of Theorem~\ref{expect} is that as long as $k
\geq c\lstar^2(T)/\theta^2$, the random operator $\Gamma/\sqrt{k}$
maps with overwhelming probability all the points in $T$ in an almost
norm preserving way.

Let us note that the method used in the proof of Theorem
\ref{expect} is called {\it generic chaining} (see \cite{Tal-book}
for the most recent survey on this subject). As we show in Section
\ref{oper_concentr}, if the set $T$ is ``very simple" one can
combine the concentration of individual variables around their
means and obtain a uniform bound using a far simpler approach.

\bigskip

The prime example for which we would like to apply
Theorem~\ref{expect} are the sets $U_m$ consisting of sparse
vectors, which are defined for $1 \le m \le n$ by
\begin{equation}
  \label{sparse}
  U_m :=\left\{x \in S^{n-1} : \left| \supp x \right| \leq m
  \right\}.
\end{equation}
We shall also consider the analogous subset of the Euclidean ball,
\begin{equation}
  \label{sparse_ball}
  {\tilde U}_m :=\left\{x \in B_2^n : \left| \supp x \right| \leq m
  \right\}.
\end{equation}
The reason for our interest in the set $U_m$ is clear: the ability to
map it in an almost norm preserving way is equivalent to the UUP. To
that end, and in light of Theorem \ref{expect}, one has to bound
$\lstar(U_m)$ in order to control $uup(\theta,\lambda)$.

The sets $U_m$ and $\tilde U_m$ have particularly simple structure:
they are the unions of the unit spheres, and unit balls, respectively,
supported on $m$-dimensional coordinate subspaces of $\R^n$.
Furthermore, for any $0 < r \le 1$,
\begin{equation}
  \label{structure_Um}
 \tilde U_{m} \cap r B_2^n = r \tilde U_{m}.
\end{equation}

It turns out that a simple geometric property of $U_m$ plays a
crucial role in the present context.

Let $T \subset \R^n$. Recall that a set $\Lambda \subset \R^n$ is an
$\eps$ cover of $T$ with respect to the Euclidean metric if
$$
T \subset \bigcup_{x \in
\Lambda} \left(x+\eps B_2^n\right),
$$
where $A+B := \{ a+b : a\in A, b \in B\}$ is the Minkowski sum of the
sets $A$ and $B$.
($\Lambda$ is often called an $\ep$-net for $T$.)
It is well-known and easy to see that if $\Lambda$ is an $\eps$ cover
of $T$ with respect to the Euclidean metric then there exists another
$\eps$ cover of $T$, say $\Lambda_1$, such that $\Lambda_1 \subset T$
and $|\Lambda_1| \le |\Lambda|$.

\medskip

The following fact is well-known and standard (see, for example,
\cite{P}, Lemma 4.10 for a part of the argument). For the convenience of
the non-specialist reader we provide a short proof.
\begin{Lemma}
   \label{lemma:ball_cover}
   Let $m \ge 1$ and $\ep >0$.  There exists an $\ep$ cover $\Lambda
   \subset B_2^m$ of $B_2^m$ with respect to the Euclidean metric
   such that $B_2^m \subset (1-\ep)^{-1} \conv \Lambda$ and $|\Lambda|
   \le (1+2/\ep)^m$.  Similarly, there exists $\Lambda' \subset
   S^{n-1}$ which is an $\ep$ cover of the sphere $S^{m-1}$ and
   $|\Lambda'| \le (1+2/\ep)^m$.
\end{Lemma}

\proof Let $\Lambda \subset B_2^m$ be a maximal subset such that
$|x-y|>\ep$ for all $x \ne y \in \Lambda$. By maximality, $\Lambda$
is an $\ep$ cover for $B_2^m$. If $x \ne y \in \Lambda$ then the two
balls $x+(\ep/2)B_2^m$ and $y+(\ep/2)B_2^m$ have disjoint interiors
and $\bigcup_{x\in \Lambda} (x+(\ep/2)B_2^m) \subset (1 +
\ep/2)B_2^m$. Comparing volumes we get the estimate for $|\Lambda|$.
Applying the same argument to the sphere $S^{m-1}$ we get a set the
desired set $\Lambda'$.  Finally, every $z \in B_2^m$ can be written
as $z = x_0 + \ep z_1$, where $x_0 \in \Lambda$ and $z_1 \in B_2^m$.
Iterating this we get that $z = x_0 + \ep x_1 + \ep^2 x_2 + \ldots
$, with $x_i \in \Lambda$, implying $B_2^m \subset (1-\ep)^{-1}
\conv\Lambda$,  as required.
\endproof

The structure of $U_m$ immediately implies similar facts as in the lemma above
for $U_m$ and $\tilde{U}_m$.
\begin{Lemma}
   \label{lemma:cover}
   There exists an absolute constant $c$ for which the following
   holds.  For every $0<\eps \leq 1/2$ and every $1 \leq m \leq n$
   there is a set $\Lambda \subset B_2^n$ which is an $\eps$ cover of
   $\tilde{U}_m$, such  that $\tilde{U}_m \subset 2 \conv \Lambda$
   and $|\Lambda| $ is at most
\begin{equation}
  \label{card_cover}
 \exp \left(c\,m \log\left(\frac{c n}{m\eps}\right)\right).
\end{equation}
Moreover, there exists an $\ep$ cover $\Lambda' \subset S^{n-1}$ of
$U_m$ with cardinality at most  (\ref{card_cover}).

Furthermore, for any $0 < r \le 1$ there exists $\bar\Lambda \subset r
B_2^n$ such that $(U_m - U_m) \cap r B_2^n \subset 2 \conv
\bar\Lambda$ and $|\bar\Lambda| $ is at most  (\ref{card_cover}).
\end{Lemma}

\proof Considering all subsets $A \subset \{1, \ldots, n\}$ with
$|A| \le m$, it is clear that the required sets $\Lambda$ and
$\Lambda'$ can be obtained as unions of the corresponding sets
supported on coordinates from $A$.  By Lemma~\ref{lemma:ball_cover}
the cardinalities of these sets are at most $(5/\eps)^m
\binom{n}{m}$.

To prove the last statement,
note that $U_m -
U_m \subset 2\tilde{U}_{2m}$, which, by (\ref{structure_Um}), implies
$$
(U_m-U_m)\cap  r B_2^n \subset r \tilde U_{2m}.
$$
By the first part of the lemma, construct a set $\Lambda\subset \tilde
U_{2m}$ such that $\tilde U_{2m} \subset 2 \conv \Lambda$ and
$|\Lambda| $  admits a suitable upper bound.
Finally, set $\bar\Lambda = r \Lambda$, completing the proof.
\endproof

\begin{Theorem}
  \label{um}
  There exist $c_1, c_2, \bar{c}, C_1 >0$ such that the following
  holds.  Let $n$, $\theta$, $X$ and $\Gamma$ be as in
  Theorem~\ref{expect}. Fix $ 1 \leq k \leq n$, let $T \subset S^{n-1}$ and assume that $T
  \subset 2 \conv \Lambda $ for some $\Lambda \subset B_2^n$  with
  $|\Lambda| \le \exp (c_1 (\theta^2/\alpha^4) k)$. Then with
  probability at least $1- \exp(- \bar{c}\, \theta^2 k/\alpha ^4)$,
  for all $x \in T$,
\begin{equation}
  \label{two_sided-3}
   1 - \theta \le \frac{|\Gamma x|^2}{  k} \le 1+\theta.
\end{equation}
Furthermore, if
$$
m \le \frac{c_2 \theta^2 k}{\alpha^4 \log(C_1 n
    \alpha^4/\theta^2 k)},
$$
then (\ref{two_sided-3}) holds for  $T = U_m$. In particular, for
every $0 < \theta <1$, and with probability at least $1- \exp(-
\bar{c}\, \theta^2 k/\alpha ^4)$, $\Gamma$
 satisfies $ uup(\theta, \lambda)$ for
$$
\lambda = \frac{ \log( n/ a' k)}{a},
$$
where both $a, a' >0$ are of the form $c \theta^2/\alpha^4$ for some
absolute constant $c$.
\end{Theorem}

The main point in the proof is that if $T \subset 2 \conv \Lambda$ for
$\Lambda \subset B_2^n$ and there is a reasonable control on the
cardinality of $\Lambda$, then $\lstar(T)$ may be bounded from above.
The rest is just a direct application of Theorem \ref{expect}.

\medskip

\proof Let $c', \bar{c} >0$ be constants from Theorem~\ref{expect}.
It is well known (see, for example, \cite{LT}) that
there exists an absolute constant $c'' >0$ such that for every
$\Lambda \subset B_2^n$,
$$
\lstar(\conv \Lambda)=\lstar (\Lambda) \leq c''\sqrt{\log(|\Lambda|)},
$$
and since $T \subset 2 \conv \Lambda$ then
$$
\lstar(T) \le 2 \lstar (\conv \Lambda) \leq
c'' \left( c_1 (\theta^2/\alpha^4) k\right)^{1/2}.
$$
Choosing $c_1 = 1/(c' c''^2)$ we conclude the proof of
(\ref{two_sided-3}) by applying Theorem~\ref{expect}.

As for the ``furthermore'' part, set $\Lambda$ to be a $1/2$ cover of
${\tilde U}_m$ provided by Lemma \ref{lemma:cover}. Then $ T = U_m
\subset {\tilde U}_m \subset 2 \conv \Lambda$. Also, by
(\ref{card_cover}) and our choice of $m$ (that includes appropriate
choices of constants $c_2$ and $C_1$), $|\Lambda |$ admits the upper
bound required in the first part of the theorem.
Finally, the last statement follows from (\ref{two_sided-3}) for
$U_m$ and Definition~\ref{uup} of the UUP.
\endproof

\section{Elementary approach}
\label{oper_concentr}

\subsection{The uniform uncertainty principle }
\label{elementary_uup}

The aim of this subsection is to obtain a positive answer to
Question~\ref{qu:CRT} using elementary methods and without resorting
to Theorem \ref{expect}. Such a proof is possible mainly because, as
already discussed in the preceding section, the geometry of the sets
$U_m$ is particularly simple.
The price one pays for the simple proof is a slightly worse dependence
on the accuracy $\theta$.

The first step in the elementary proof is obtaining an analog of
Theorem \ref{expect}, where the complexity measure $\lstar(T)$ is
replaced by estimates on covering numbers.

Consider a set of random $k \times n$ matrices $\tilde\Gamma$
satisfying two conditions.
First,
\begin{equation}
  \label{gen-1}
    \E |\tilde\Gamma x|^2 = 1
\qquad \mbox{\rm   for all }x \in S^{n-1},
\end{equation}
that is, on average, $\tilde \Gamma$ preserves the norm of each
individual $x$.

The second condition asserts the concentration of the random
variable $|\tilde\Gamma x|^2$ around its expectation: there exists
an absolute constant $c_0$ such that for every $x \in \R^n$ we have
\begin{equation}
  \label{gen-2}
  \Pp \left(\bigl | |\tilde\Gamma x|^2 -|x|^2 \bigr| \ge  t |x|^2 \right)
\le  e^{- c_0 t^2 k} \qquad \mbox{\rm for all } 0 < t \le  1.
\end{equation}

Let us note that (multiples of) subgaussian matrices considered in
Section~\ref{subgaussian} satisfy (\ref{gen-1}) and (\ref{gen-2}).
Indeed, let $(X_i)_{i=1}^k$ be independent copies of an isotropic
$\psi_2$ vector with a constant $\alpha$ and set
$$
\tilde \Gamma = \frac{1}{\sqrt{k}}
\sum_{i=1}^k \inr{X_i, \cdot} e_i,
$$
where $(e_i)_{i=1}^n$ are the standard unit vectors in $\R^n$. By
the isotropicity assumption, $\E|\tilde \Gamma x|^2=1$ for every $x
\in S^{n-1}$. Moreover, by fixing $x \in S^{n-1}$ and applying
Bernstein's inequality (see, e.g. \cite{LT,VW}) to the average of
$k$ independent copies of the random variable $\inr{X,x}^2$, it
is evident
that for every $t>0$,
$$
\Pp \left( \left| \frac{1}{k}\sum_{i=1}^k \inr{X_i,x}^2 -1 \right|
> t \right) \leq
2\exp\left(-ck\min\left\{\frac{t^2}{\alpha^4},
              \frac{t}{\alpha^2}\right\}\right),
$$
where $c$ is an absolute constant. Since $\alpha \geq 1$, $\tilde
\Gamma$ satisfies \eqref{gen-2} for $c_0=c/\alpha^4$.

\bigskip
Let us formulate the elementary version of Theorem \ref{expect}.
\begin{Theorem}
  \label{elementary}
  Consider a set of random $k \times n$ matrices $\tilde\Gamma$
  satisfying (\ref{gen-1}) and (\ref{gen-2}).  Let $T \subset S^{n-1}$
  and $0 < \theta < 1$, and assume the following:
\begin{description}
\item{(i)} There exists $\Lambda' \subset S^{n-1} $ which is a
  $\theta/5$-cover of $T$ and satisfies $|\Lambda'| \le \exp( c_0
  \theta^2 k/50)$.
\item{(ii)} There exists $\Lambda \subset (\theta/5) B_2^n$ such that
  $(T-T)\cap (\theta/5) B_2^n \subset 2 \conv \Lambda$ and $|\Lambda|
  \le \exp( c_0 k/2)$.
\end{description}
Then with probability at least $1 - 2 \exp( -c_0 \theta^2 k/50)$,
for all $x \in T$,
\begin{equation}
  \label{two_sided-4}
   1 - \theta \le |\tilde\Gamma x|^2   \le 1+\theta.
\end{equation}
\end{Theorem}

\begin{Remark} {\rm There is nothing special
in  the constant $2$ in front of $\conv \Lambda$ in (ii), and
    it could be replaced by any constant strictly larger than $1$.}
\end{Remark}

The idea behind the proof of Theorem \ref{elementary} is to show
that $\tilde \Gamma$ acts on $\Lambda '$ in an almost norm
preserving way. This is the case because the degree of concentration
of each variable $|\tilde \Gamma x|^2$ around its mean defeats the
cardinality of $\Lambda '$. Then one shows that $\tilde \Gamma
(\conv \Lambda )$ is contained in a small ball - thanks to a similar
argument.

\medskip

\proof Set $\ep = \theta/5$ and consider the set of $\tilde\Gamma$ on
which
\begin{equation}
  \label{inf}
\left||\tilde\Gamma x_0 |-1 \right| \le \left||\tilde\Gamma x_0 |^2 -
1 \right| \le \ep \qquad \mbox{\rm for \ all \ } x_0 \in \Lambda',
\end{equation}
and
\begin{equation}
  \label{sup}
|\tilde\Gamma z | \le 2 |z|
\qquad \mbox{\rm for \ all \ } z \in \Lambda.
\end{equation}
Note  that this set has probability larger than or equal to $ 1 -
\exp(-c_0 \ep^2 k/2) -  \exp(-c_0 k/2)
\ge 1 - 2 \exp(-c_0 \ep^2 k/2)$.

\smallskip

Let $x \in T$ and consider $x_0 \in \Lambda'$ such that $| x -
x_0| \le \ep$. Then
$$
|\tilde\Gamma x_0 | - |\tilde\Gamma (x-x_0)| \le |\tilde\Gamma x|
\le |\tilde\Gamma x_0 | + |\tilde\Gamma (x-x_0)|.
$$
Since
 $x - x_0 \in (T-T)\cap \ep  B_2^n $, then by the definition
of $\Lambda$
and (\ref{sup}) it follows that
\begin{equation}
  \label{sup-2}
  |\tilde\Gamma (x-x_0) |\le 2 \sup_{z \in \Lambda} |\tilde\Gamma z|
  \le 4 \ep.
\end{equation}
Combining this with (\ref{inf}) implies that
$1- 5 \ep \le |\tilde\Gamma x | \le 1 + 5 \ep$, completing the
proof, by the definition of $\ep$.
\endproof

We are now ready for an elementary solution to Question \ref{qu:CRT},
contained in the following corollary.

\begin{Corollary}
  \label{quest-1}
  Let $\tilde \Gamma$ satisfy \eqref{gen-1} and \eqref{gen-2}. Then,
  there are constants $c_1, c_1'$ and $c_2$ depending only on $c_0$
  from \eqref{gen-2} for which the following holds. For every
  $0<\theta<1$, with probability at least $1- 2\exp(-c_2\theta^2k)$,
  $\tilde \Gamma$ satisfies $uup(\theta,\lambda)$ for
$$
\lambda=\frac{c_1 \log\left(c_1' n/k\theta^3\right)}{\theta^2}.
$$
In particular, there is are absolute constants $c_1$, $c_2$ and
$c_3$ for which the following holds. If $X$ is an isotropic,
$\psi_2$ vector with constant $\alpha$ then with probability at
least $1-\exp(-c_1\theta^2k/\alpha^4)$, the operator $\tilde
\Gamma  =\frac{1}{\sqrt{k}}\sum_{i=1}^k \inr{X_i,\cdot}e_i$
satisfies $uup(\theta,\lambda)$ for
$$
\lambda=\frac{c_2\alpha^4}{\theta^2
\log\left(c_3n\alpha^4/k\theta^3\right)}.
$$
\end{Corollary}


\medskip

\proof The main part of the proof is to show that there exists $c' >0$
such that, given $0<\theta<1$, if $m$ and $k$ satisfy
\begin{equation}
  \label{k_and_m}
k \geq \frac{c'm}{\theta^2}\log \left(\frac{c'n}{m \theta}\right),
\end{equation}
then (\ref{two_sided-4}) holds, that is, $\tilde \Gamma$ acts on
$U_m$ in an almost norm preserving way.  To that end we need to
exhibit the sets $\Lambda$ and $\Lambda'$.  For the latter set, one
can choose $c'$ in such a way that the set $\Lambda'$ constructed in
the moreover part of Lemma~\ref{lemma:cover} for $\ep = \theta/5$
satisfies the required condition (i) for $T = U_m$.  For the former
set, apply the third part of Lemma~\ref{lemma:cover} with $r =
\theta/5$ to get $\bar\Lambda$; adjusting the choice of $c'$ in
(\ref{k_and_m}), $\bar\Lambda$ satisfies (ii).


Now the conclusion follows from (\ref{k_and_m}) by a straightforward
computation.
\endproof

\medskip

\begin{Remark} {\rm Note that the price for using the elementary
approach in the case of $U_m$ - and thus for Question \ref{qu:CRT}
is not very high - a slightly worse power of $\theta$ in the
logarithm. However, there are many cases of sets $T \subset S^{n-1}$
in which this elementary approach would not be enough to show that
$\tilde \Gamma$ acts in an almost norm preserving way on $T$.}
\end{Remark}

\subsection{The approximate reconstruction problem}
\label{approximate}

Next, we show how the elementary approach can be used to solve the
approximate reconstruction problem in several cases that have been
considered in \cite{CT1, CT2, D, MPT1, MPT2, BDDW}, among others.  Let
us recall the formulation of this problem.

\begin{Question}
\label{aprox-question}
Suppose that $\tilde T \subset \R^n$ and fix $t_0 \in \tilde T$. Let
$\Gamma$ be a $k \times n$ random matrix, and suppose that one is
given the data vector $\Gamma t_0$, that is, the set of linear
measurements $\left(\inr{X_i,t_0}\right)_{i=1}^k$. Is it possible to
find (with high probability) some $x \in \R^n$, such that $|x-t_0|$ is
small?
\end{Question}

In \cite{CT1} this problem has been studied by using the UUP and for
particular sets -- $B_1^n$, the unit ball in $\ell_1^n$ and
$B_{p,\infty}^n$ for $0<p<1$, the unit balls in weak $\ell_p$ spaces.
In \cite{MPT1,MPT2}, a geometric approach was introduced which solved
this problem for an arbitrary symmetric quasi-convex subset of $\R^n$.
(Recall that a (centrally) symmetric set $\tilde T$ is quasi-convex
with constant $a \ge 1$, if $\tilde T + \tilde T \subset 2 a \tilde T$
and $\tilde T$ is star-shaped, i.e., $s \tilde T \subset \tilde T$ for
$0 < s < 1$.)

The geometric idea at the heart of \cite{MPT1,MPT2} is essentially the
following: let $\tilde T \subset \R^n$ and suppose that one can find
$\eps_k$ and show that with high probability,
$$
{\rm diam}\left(\ker(\Gamma)\cap \tilde T\right) \leq \eps_k.
$$
Since $\tilde T$ is quasi-convex,
then $\tilde T-\tilde T \subset 2a \tilde T$, for some $a \ge 1$.
Hence, ${\rm diam}\left(\ker(\Gamma)\cap (\tilde T- \tilde T)\right)
\leq 2a \eps_k$.
In particular, with high probability, if $x $ is in $\tilde T$ and it
satisfies $\Gamma x =\Gamma t_0$ then $|x-t_0| \leq 2 a \eps_k$, as
required.

In other words, the approximate reconstruction problem is reduced to
finding an upper estimate on the diameter of the intersection of the
kernel of $\Gamma$ with $\tilde T$ that holds with high probability.
This parameter has been studied in asymptotic geometry and in
approximation theory for certain notions of randomness, and is the
random $k$-th {\it Gelfand number} of $T$ associated with the random
matrix $\Gamma$.

\medskip

Theorem~\ref{elementary} provides a method for estimating the
diameter of $\ker(\Gamma) \cap \tilde T$ in the following way.  For
$\tilde T \subset \R^n$ star-shaped let $T_\rho = \tilde T \cap \rho
S^{n-1}$. Then if (a multiple of) $\Gamma$ acts on $T_\rho$ in an
almost norm preserving way, then $\ker(\Gamma) \cap \tilde T \subset
\rho B_2^n$, and thus ${\rm diam}\left(\ker(\Gamma)\cap \tilde
  T\right) \leq \rho$.  Indeed, if not, there would be a point $t \in
\tilde T$ of norm greater than $\rho$ which is mapped to $0$. Hence,
$\rho t/|t| \in T$ will also be mapped to $0$, which contradicts the
fact that (a multiple of) $\Gamma$ is almost norm preserving on
$\tilde T$.

This proves  the following Corollary.

\begin{Corollary}
  \label{diam}
  Let $\tilde\Gamma$ be as in Theorem~\ref{elementary}.  Let $\tilde T
  \subset \R^n$ be star-shaped. Let $T = \rho^{-1}(\tilde T \cap \rho S^{n-1})$ and
  assume that $T $ satisfies the hypothesis of
  Theorem~\ref{elementary} for some $0 < \theta< 1$ (say, $\theta =
  1/2$).  Then ${\rm diam}\left(\ker(\tilde \Gamma)\cap \tilde
    T\right) \leq \rho$, with probability at least $1 - 2 \exp (-ck)$,
  where $c >0$ is an absolute constant.
\end{Corollary}

\bigskip

To illustrate this corollary, we consider examples of $\tilde T$:
the unit ball in $\ell_1^n$, denoted by $B_1^n$, and the unit balls
in $\ell_p^n$ and the weak-$\ell_p^n$ spaces $\ell_{p,\infty}^n$ for
$0<p<1$, denoted by $B_p^n$ and $B_{p,\infty}^n$, respectively.
Recall that $B_{p,\infty}^n$ is the set of all $x = (x_i)_{i=1}^n
\in \R^n$ such that the cardinality $|\{i: |x_i|\ge s\}|\le s^{-p}$
for all $ s>0$. Note that $B_p^n \subset B_{p, \infty}^n$ so we can
restrict ourselves to considering the balls $B_{p, \infty}^n$ only.

We will require two lemmas.  The first lemma comes from \cite{MPT2}
and it combines a reformulation of Lemma 3.2 and (3.1) from that
article.

\begin{Lemma}
  \label{convex_hull}
Let $0<p<1$, $1 \leq m \leq n$ and set
$ r =(1/p -1 ) m^{1/p- 1/2}$. Then, for every $x \in \R^n$,
$$
\sup_{z \in r B_{p,\infty}^n \cap  B_2^n} \inr{x,z} \leq 2
\left(\sum_{i=1}^m {x_i^*}^2\right)^{1/2},
$$
where $(x_i^*)_{i=1}^n$ is a non-increasing rearrangement of
$(|x_i|)_{i=1}^n$. Equivalently,
\begin{equation}
  \label{convex_B_p}
r B_{p,\infty}^n \cap B_2^n \subset 2 \conv \tilde U_m.
\end{equation}
Furthermore,
\begin{equation}
  \label{convex_B_1}
\sqrt m  B_1^n \cap B_2^n  \subset 2 \conv \tilde U_m.
\end{equation}
\end{Lemma}

The second lemma shows that $m^{1/p-1/2}B_{p,\infty}^n\cap S^{n-1}$
is well approximated by vectors on the sphere with a relatively
short support.

\begin{Lemma}
  \label{net}
  Let $0<p<2$ and $\delta >0$, set $\ep = 2(2/p-1)^{-1/2}
  \delta^{1/p-1/2}$. Then $ U_{\lceil{m/\delta} \rceil}$ is
  an $\ep$-cover of
  $m^{1/p-1/2}B_{p,\infty}^n\cap S^{n-1}$ with respect to the
  Euclidean metric.
\end{Lemma}

\proof
Let $x \in m^{1/p-1/2}B_{p,\infty}^n\cap S^{n-1}$ and
assume without loss of generality that $x_1\ge x_2\ge \ldots \ge x_n
\ge 0$. Define $z'$ by $z_i = x_i$ for $1 \le i \le \lceil{m/\delta}
\rceil$ and $z'_i = 0$, otherwise. Then
$$
|x-z'|^2 = \sum_{i>m/\delta} |x_i|^2
   \le m^{2/p-1}\sum_{i>{m/\delta}} {1/ i^{2/p}}
\le (2/p -1 )^{-1}\, \delta^{2/p-1}.
$$
Thus $ 1 \ge |z'|\ge 1 - (2/p -1 )^{-1/2}\, \delta^{1/p-1/2}$.  Put $z
= {z'}/{|z'|}$. Then $z \in U_{\lceil{m/\delta} \rceil}$ and
$$
|z - z'| =  1 - |z'| \le (2/p -1 )^{-1/2}\, \delta^{1/p-1/2}.
$$
By the triangle inequality $|x-z|\le \ep$,
completing the proof.
\endproof

\bigskip



Let $0 < p < 1$.  Fix $1 \le m \le n$, set $\tilde T :=
m^{1/p-1/2}B_{p,\infty}^n$ and $T := \tilde T \cap S^{n-1}$. We
shall show that for appropriately chosen $m$, $T $ satisfies the
hypothesis of Theorem~\ref{elementary} for $\theta = 1/2$. To that
end, we need to show that the complexity of the set $T$ as captured
by the sets $\Lambda$ and $\Lambda'$ is small.

First note, to simplify the calculations a little, that by
Lemma~\ref{net}, for $\delta >0$ the set $ U_{\lceil{m/\delta}
  \rceil}$ is an $\ep$ cover for $T$, where $\ep = 2 \sqrt
\delta$. (That is, the dependence of $\ep$ on $\delta$ is universal
in the range of $p$ considered here.)  Use this fact for $\delta =
1/40^2$ and combine it with the ``moreover part'' of
Lemma~\ref{lemma:cover} (for $\ep = 1/20$) which provides us with a
set $\Lambda' \subset S^{n-1}$ which is $1/20$ cover of $
U_{\lceil{m/\delta}  \rceil}$. Hence, by the triangle inequality,
$\Lambda'$ is $1/20 + 1/20 = 1/10$ cover of $T$. Moreover, by
(\ref{card_cover}), $|\Lambda'| \le \exp\left( c_1 m \log
  (c_1 n/m)\right)$, where $c_1>0 $ is an absolute constant.

It is easy to check that that $B_{p, \infty}$ is quasi-convex with
constant $2^{1/p}$ and therefore
$$
(T-T)\cap {{1}\over{10}} B_2^n \subset
\left(2^{1+1/p} \tilde T \cap 2B_2^n\right)\cap  {1\over 10} B_2^n =
2^{1+1/p} \tilde T \cap  {1\over 10} B_2^n = {1\over 10} A,
$$
where
$$
A:= \left(10\,\cdot\, 2^{1+1/p}\right) \, \tilde T \cap B_2^n =
\left(10\,\cdot\, 2^{1+1/p} \, m^{1/p-1/2}\right) B_{p,\infty}^n \cap
B_2^n.
$$
Set $m_1 = \max\left(c_p'm, m\right) $ where ${c_p'}^{1/p - 1/2}= (1/p
-1)^{-1} 20 \cdot 2^{1/p}$, so that
$$
10\,\cdot\, 2^{1+1/p} \, m^{1/p-1/2} \le
(1/p-1) m_1^{1/p - 1/2}.
$$
Then, by (\ref{convex_B_p}), $A \subset 2 \conv \tilde U_{m_1}$. By
the first part of Lemma~\ref{lemma:cover} there is a subset
$\Lambda_1\subset B_2^n$ such that $ \tilde U_{m_1} \subset 2 \conv
\Lambda_1$ and $|\Lambda_1|\le \exp\left( c_1' m_1  \log (c_1'
n/m_1)\right)$, where $c_1' >0$ is an absolute constant. Letting
$\Lambda = {1\over 10} \Lambda_1$ yields $(T-T)\cap {{1}\over{10}}
B_2^n \subset 4\conv \Lambda$ and $|\Lambda|\le \exp\left( c_p'' m
\log (c_1' n/m)\right)$, where  $c_p'' \ge 1 $ depends on $p$ only.
(The precise form of $c_p''$  can be easily calculated from the form
of $c_p'$ but we shall not do it here.)

Considering the upper bounds for $\Lambda'$ and $\Lambda$ yields the
existence of $c_p \ge 1$, depending on $p$ only, and of an absolute
constant $c_1' >0$
such that whenever $k$ satisfies
\begin{equation}
  \label{k-m-two}
k \geq c_p m\log \left(\frac{c_1' n}{m}\right),
\end{equation}
then $T = m^{1/p-1/2}B_{p,\infty}^n \cap S^{n-1}$ satisfies
assumptions (i) and (ii) of Theorem~\ref{elementary}.  Therefore, by
Corollary~\ref{diam}, ${\rm diam}\left(\ker(\tilde \Gamma)\cap
  B_{p,\infty}^n \right) \leq m^{1/2-1/p}$, with high probability.

For $p=1$, an analogous result holds for $B_1^n$:
if $k$ and $m$ satisfy (\ref{k-m-two}) (with  $c_p$ replaced by a
certain absolute constant)  then
 ${\rm diam}\left(\ker(\tilde \Gamma)\cap
  B_{1}^n \right) \leq m^{-1/2}$, with high probability.

 A straightforward
calculation then leads to the following estimates
for the diameters of the intersection of $\ker(\tilde \Gamma)$
with the balls $B_p^n$ and $B_{p, \infty}^n$ (for $0 < p <1$) and of
$B_1^n$.
\begin{Corollary}
  \label{diam_weak-lp}
  Let $\tilde\Gamma$ be as in Theorem~\ref{elementary}.  Let $0 < p
  <1$. There exist a constant $c_p$  depending  only on $p$, a
  constant $c$  depending on $c_0$ and an absolute constant $c_1$,
  such that, with probability at least $1-\exp(-ck)$,
$$
\diam\left(\ker(\tilde\Gamma) \cap B_{p}^n  \right) \leq
\diam\left(\ker(\tilde\Gamma) \cap B_{p,\infty}^n  \right) \leq
c_p\left(\frac{\log(c_1n/k)}{k}\right)^{1/p - 1/2}.
$$
In particular, of $t_0 \in B_{p,\infty}^n$ and
$\tilde{\Gamma}x=\tilde{\Gamma} t_0$ then with high probability,
$$
|x-t_0| \leq c'_p\left(\frac{\log(c_1n/k)}{k}\right)^{1/p - 1/2}.
$$
For $p=1$, an analogous result holds for the ball $B_1^n$ replacing
$B_{p, \infty}^n$ and $c_p$ and $c'_p$ being replaced by an
absolute constant.
\end{Corollary}



\footnotesize
{
}

\noindent {\bf S. Mendelson} {\footnotesize Centre for Mathematics
and its Applications, The
Australian National University, Canberra, ACT 0200,
Australia \\} {\small\tt%
shahar.mendelson@anu.edu.au}  \\ [.05cm]

\noindent {\bf A. Pajor }{\footnotesize Laboratoire d'Analyse et
Math\'ematiques Appliqu\'ees, Universit\'e de Marne-la-Vall\'ee, 5
boulevard Descartes, Champs sur Marne, 77454 Marne-la-Vallee,
Cedex 2, France \\ }
{\small\tt%
    alain.pajor@univ-mlv.fr}\\ [.05cm]

\noindent {\bf N. Tomczak-Jaegermann} {\footnotesize
Department of Mathematical and Statistical Sciences,\\
University of Alberta,
Edmonton, Alberta, Canada T6G 2G1\\ }
{\small\tt%
  nicole@ellpspace.math.ualberta.ca}

\end{document}